\documentclass[12pt,reqno]{amsart}

\usepackage{amssymb}
\usepackage[arrow]{xypic}

\DeclareMathOperator{\chr}{char}\DeclareMathOperator{\decc}{dec}
\DeclareMathOperator{\Fix}{Fix} \DeclareMathOperator{\Gal}{Gal}
 
\DeclareMathOperator{\Ker}{ker} \DeclareMathOperator{\proj}{proj}
\DeclareMathOperator{\res}{res} \DeclareMathOperator{\tra}{tra}

\newcommand{\dec}{{\decc}}
\newcommand{\Ec}{\mathcal{E}}
\newcommand{\F}{\mathbb{F}}
\newcommand{\Fp}{\F_p}
\newcommand{\Ft}{\F_2}
\newcommand{\G}{{\Gamma}}
\newcommand{\GG}{\mathcal{G}}
\renewcommand{\H}{\mathcal{H}}

\newcommand{\Ic}{\mathcal{I}}
\newcommand{\Jc}{\mathcal{J}}
\newcommand{\Kc}{\mathcal{K}}
\newcommand{\N}{{\Delta}}
\newcommand{\Nb}{\mathbb{N}}

\newcommand{\Z}{\mathbb{Z}}

\begin{document}

\title[Detecting Pro-$p$-groups]{Detecting Pro-$p$-groups that
are not Absolute Galois Groups}

\author[Dave Benson]{Dave Benson}
\address{Department of Mathematical Sciences, University of
Aberdeen, Meston Building, King's College, Aberdeen AB24 3UE,
Scotland UK}
\email{bensondj@maths.abdn.ac.uk}

\author[Nicole Lemire]{Nicole Lemire}
\author[J\'{a}n Min\'{a}\v{c}]{J\'an Min\'a\v{c}}
\address{Department of Mathematics, Middlesex College, \ University
of Western Ontario, London, Ontario \ N6A 5B7 CANADA}
\email{nlemire@uwo.ca}
\email{minac@uwo.ca}

\author[John Swallow]{John Swallow}
\address{Department of Mathematics, Davidson College, Box 7046,
Davidson, North Carolina \ 28035-7046 USA}
\email{joswallow@davidson.edu}

\date{October 20, 2006}

\thanks{Nicole Lemire was supported in part by Natural Sciences
and Engineering Research Council of Canada grant R3276A01. J\'an
Min\'a\v{c} was supported in part by Natural Sciences and
Engineering Research Council of Canada grant R0370A01.  John Swallow
was supported in part by National Science Foundation grant
DMS-0600122.}

\maketitle

\def\thepart{A}

\newtheorem{theorem}{Theorem}[part]
\newtheorem{proposition}{Proposition}[part]
\newtheorem{corollary}[theorem]{Corollary}
\newtheorem*{corollary*}{Corollary}
\newtheorem{lemma}{Lemma}[part]

\theoremstyle{definition}
\newtheorem*{remark*}{Remark}
\newtheorem*{example*}{Example}

\parskip 9pt plus 3pt

Let $p$ be a prime.  It is a fundamental problem to classify the
absolute Galois groups $G_F$ of fields $F$ containing a primitive
$p$th root of unity $\xi_p$.  In this paper we present several
constraints on such $G_F$, using restrictions on the cohomology of
index $p$ normal subgroups from \cite{LMS}. In section~\ref{se:t} we
classify all maximal $p$-elementary abelian-by-order $p$ quotients
of these $G_F$.  In the case $p>2$, each such quotient contains a
unique closed index $p$ elementary abelian subgroup. This seems to
be the first case in which one can completely classify nontrivial
quotients of absolute Galois groups by characteristic subgroups of
normal subgroups. In section~\ref{se:asb} we derive analogues of
theorems of Artin-Schreier and Becker for order $p$ elements of
certain small quotients of $G_F$. Finally, in section~\ref{se:fam}
we construct a new family of pro-$p$-groups which are not absolute
Galois groups over any field $F$.

As a consequence of our results, we prove the following limitations
on relator shapes of pro-$p$ absolute Galois groups. For elements
$\sigma$ and $\tau$ of a group $\G$, let ${}^0[\sigma,\tau]=\tau$,
${}^1[\sigma,\tau]=\sigma\tau\sigma^{-1}\tau^{-1}$, and
${}^n[\sigma,\tau] = [\sigma, {}^{n-1}[\sigma,\tau]]$ for $n\ge 2$.
Similarly, for subsets $\G_1$ and $\G_2$ of $\G$, let
${}^n[\G_1,\G_2]$ denote the closed subgroup generated by all
elements of the form ${}^n[\gamma_1, \gamma_2]$ for $\gamma_i\in
\G_i$.
\begin{theorem}\label{th:1}
    Let $p$ be an odd prime, $\G$ a pro-$p$-group with maximal
    closed subgroup $\N$, and $\sigma\in \G\setminus \N$.
    \begin{enumerate}
        \item Suppose that for some $\tau\in \N$ and some $2\le e\le
        p-2$
        \begin{equation*}
            {}^e[\sigma,\tau]\not\in {}^{p-1}[\sigma,\N]\Phi(\N)
            \quad\text{and}\quad
            {}^{e+1}[\sigma,\tau] \in \Phi(\N).
        \end{equation*}
        Then $\G$ is not an absolute Galois group.

        \

        \noindent
        Moreover, if $\G$ contains a normal subgroup $\Lambda\subset
        \N$ such that $\G/\Lambda\simeq \Z/p^2\Z$, we may take $1\le
        e\le p-2$.
        \item Suppose that for some $\tau_1, \tau_2\in \N$,
        \begin{equation*}
            [\sigma,\tau_i]\not\in {}^{p-1}[\sigma,\N]\Phi(\N)
            , \ {}^2[\sigma,\tau_i] \in \Phi(\N), \quad i=1,
            2
        \end{equation*}
        \begin{equation*}
            \langle [\sigma,\tau_1] \rangle \Phi(\N) \neq
            \langle [\sigma,\tau_2] \rangle \Phi(\N)
        \end{equation*}
        Then $\G$ is not an absolute Galois group.
        \item Suppose that
        \begin{equation*}
            \sigma^p \in {}^{2}[\sigma,\N]\Phi(\N).
        \end{equation*}
        Then $\G$ is not an absolute Galois group.
    \end{enumerate}
\end{theorem}
\noindent Here $\Phi(\N)=\N^p[\N,\N]$ denotes the Frattini subgroup
of $\N$.

Furthermore, pro-$p$-groups with single relations similar to those
of Demu\v{s}kin groups for odd primes cannot be absolute Galois
groups.
\begin{corollary*}\label{co:1}
    Let $p$ be an odd prime and $\G$ a pro-$p$-group minimally
    generated by $\{\sigma_1,\sigma_2\} \cup \{\sigma_i\}_{i\in
    \Ic}$ subject to a single relation
    \begin{equation*}
        \sigma_1^q\cdot {}^f[\sigma_1,\sigma_2]\cdot
        \prod_{(i,j)\in \Jc} [\sigma_i,\sigma_j]\cdot
        \prod_{k \in \Kc} [\sigma_1^p,\sigma_k]
    \end{equation*}
    for some $2\le f\le p-1$, $q\in \Nb\cup \{0\}$ with $q = 0 \bmod
    (p^2)$, a finite ordered set of pairs $\Jc\subset \Ic\times
    \Ic$, and a finite ordered subset $\Kc$ of $\Ic$.  Then $\G$ is
    not an absolute Galois group.
\end{corollary*}

The results in \cite{LMS} may be used to establish further new
results on possible $V$-groups of fields and metabelian quotients of
absolute Galois pro-$p$-groups. (For the definition of the $V$-group
$V_F$ of a field $F$, see section~\ref{se:asb}.) Moreover, some of
the results here hold in a greater generality than their
formulations here.  For instance, the examples here of pro-$p$
non-absolute Galois groups are also examples of groups which are not
maximal pro-$p$-quotients of absolute Galois groups, by \cite[\S
6]{LMS}. Furthermore, pro-$p$-groups which are not absolute Galois
groups are not $p$-Sylow subgroups of absolute Galois groups.  We
plan a systematic study of these concerns in \cite{BeMS}.

We observe that because this paper is concerned only with degree $1$
and degree $2$ cohomology, the results cited from \cite{LMS} rely
only on the Merkurjev-Suslin Theorem \cite[Theorem 11.5]{MeSu} and
not the full Bloch-Kato Conjecture. Furthermore, we note that while
this paper is self-contained, an extended version is available
\cite{longapp}.

\section{$T$-groups}\label{se:t}

A \emph{$T$-group} is a nontrivial pro-$p$-group $T$ with a maximal
closed subgroup $N$ that is abelian of exponent dividing $p$. Then
$N$ is a normal subgroup, and the factor group $T/N$ acts naturally
on $N$: choose a lift $t\in T$ and act via $n\mapsto tnt^{-1}$.
Given any profinite group $\G$ with a closed normal subgroup $\N$ of
index $p$, the factor group $\G/\N^p[\N,\N]$ is a $T$-group. Now
suppose that $E/F$ is a cyclic field extension of degree $p$. We
define the \emph{$T$-group of $E/F$} to be $T_{E/F} :=
G_F/G_E^p[G_E,G_E]$.  In this section we classify those $T$-groups
realizable as $T_{E/F}$ for fields $F$ either with $\chr F=p$ or
$\xi_p\in F$.

We develop a complete set of invariants $t_i$, $i=1, 2, \dots, p$,
and $u$ of $T$-groups as follows. For a pro-$p$-group $\G$, denote
by $Z(\G)$ its center, $Z(\G)[p]$ the elements of $Z(\G)$ of order
dividing $p$, $\G^{(n)}$ the $n$th group in the $p$-central series
of $\G$, and $\G_{(n)}$ the $n$th group in the central series of
$\G$. For a $T$-group $T$ we define
\begin{align*}
    t_1 &= \dim_{\Fp} H^1\left(\frac{Z(T)[p]}{Z(T)\cap T_{(2)}},
    \Fp\right) \\ t_i &= \dim_{\Fp} H^1 \left(\frac{{Z(T)\cap
    T_{(i)}}}{Z(T)\cap T_{(i+1)}},\Fp\right), \quad 2\le i\le p\\
    u &= \max\{i: 1\le i\le p, \ T^p\subset T_{(i)}\}.
\end{align*}

\begin{proposition}\label{pr:tc}
    For arbitrary cardinalities $t_i$, $i=1, 2, \dots, p$, and
    $u$ with $1\le u\le p$, the following are equivalent:
    \begin{enumerate}
        \item\label{it:tc1} The $t_i$ and $u$ are invariants of some
        $T$-group
        \item\label{it:tc2}
        \begin{enumerate}
            \item\label{it:tc2a} if $u<p$ then $t_u\ge 1$, and
            \item\label{it:tc2b} if $u=p$ and $t_i=0$ for all $2\le
            i\le p$, then $t_1\ge 1$.
        \end{enumerate}
    \end{enumerate}
    Moreover, $T$-groups are uniquely determined up to isomorphism
    by these invariants.
\end{proposition}

\begin{theorem}\label{th:t} For $p$ an odd prime, the
    following are equivalent:
    \begin{enumerate}
        \item $T$ is a $T$-group with invariants $t_i$ and $u$
        satisfying
        \begin{enumerate}
            \item $u\in \{1,2\}$,
            \item $t_2=u-1$, and
            \item $t_i=0$ for $3\le i<p$.
        \end{enumerate}
        \item $T\simeq T_{E/F}$ for some cyclic extension $E$ of degree
        $p$ of a field $F$ such that either $\chr F=p$ or $\xi_p\in F$.
    \end{enumerate}

    Now suppose $p=2$.  Then each $T$-group is isomorphic to
    $T_{E/F}$ for some cyclic extension $E/F$ of degree $2$.
\end{theorem}

Let $G$ be a cyclic group of order $p$ and $M_i$, $i=1, 2, \dots,
p$, denote the unique cyclic $\Fp G$-module of dimension $i$.  Since
the $\Fp G$-modules we consider will be multiplicative groups, we
usually write the action of $G$ exponentially. For a set $\Ic$, let
$M_i^{\Ic}$ denote the product of $\vert\Ic\vert$ copies of $M_i$
endowed with the product topology.  We use the word \emph{duality}
exclusively to refer to Pontrjagin duality, between compact and
discrete abelian groups and more generally between compact and
discrete $\Fp G$-modules.  We denote the dual of $\Gamma$ by
$\Gamma^\vee$.

\begin{lemma}\label{le:t}
    Suppose $T$ is a $T$-group with invariants $t_i$, $i=1, 2,
    \dots, p$ and $u$, and $N$ is a maximal closed subgroup of $T$
    that is abelian of exponent dividing $p$.  Let $\sigma\in
    T\setminus N$ and set $\bar\sigma$ to be the image of $\sigma$
    in $G:=T/N$.  Then $N$ is a topological $\Fp G$-module and we
    have
    \begin{enumerate}
        \item\label{it:lt2} For any $\Fp G$-submodule $M$ of $N$ and
        $i\ge 0$, ${}^i[T,M]=M^{(\bar\sigma-1)^i}$.
        \item\label{it:lt3} For $i\ge 2$, $T_{(i)} =
        N^{(\bar\sigma-1)^{i-1}}$.
        \item\label{it:lt4} There exist sets $\Ic_i$, $i=1, 2,
        \dots, p$, such that $N$ decomposes into indecomposable $\Fp
        G$-modules as $N = M_1^{\Ic_1}\times M_2^{\Ic_2} \times
        \cdots \times M_p^{\Ic_p}$, endowed with the product
        topology. Moreover, for $i\ge 2$, $t_i=\vert \Ic_i\vert$,
        and
        \begin{equation*}
            t_1 = \begin{cases}
                1+\vert \Ic_1\vert, &T \text{\ is abelian\ of\
                exponent\ } p\\
                \vert \Ic_1\vert, &\text{otherwise}.
            \end{cases}
        \end{equation*}
        \item\label{it:lt5} $T^p = \langle \sigma^p \rangle \cdot
        T_{(p)}$.
        \item\label{it:lt6} If $u<p$ then $u$ is the minimal $v$ with
        $1\le v\le p-1$ such that there is a commutative diagram of
        pro-$p$-groups
        \begin{equation}\label{eq:cd}
            \vcenter{\xymatrix{ 1 \ar[r]& N \ar[r] \ar@{>>}[d] & T
            \ar[r] \ar@{>>}[d] & G \ar[r] \ar[d]^{=} & 1 \\ 1 \ar[r]
            & M_v \ar[r] & H \ar[r] & G \ar[r] & 1}}
        \end{equation}
        with a lift of $\bar\sigma$ in $H$ of order $p^2$. If $u=p$
        then no such diagram exists for $1\le v\le p-1$.
    \end{enumerate}
\end{lemma}

\begin{proof}
    (\ref{it:lt2}).  Suppose that $\tau\in T$ is arbitrary, and
    write $\tau= n\sigma^{i}$ for $n\in N$ and $i \in \Nb \cup
    \{0\}$.  The action of $T$ factors through $G$.  Hence
    $[\tau,M] = M^{(\bar\sigma^i-1)}$ and so $[T,M]=
    M^{(\bar\sigma-1)}$.  The result follows by induction.
   (Moreover, we observe that if $p$ does not divide $i$, then
    $[\tau,M] =[M,\tau]=[T,M]=[M,T]$.)

    (\ref{it:lt3}).  Observe that $[T,T]=[T,N]$.  Then use
    \eqref{it:lt2}.

    (\ref{it:lt4}).  Because $\Fp G$ is an Artinian principal ideal
    ring, every $\Fp G$-module $U$ decomposes into a direct sum of
    cyclic $\Fp G$-modules.  Every cyclic $\Fp G$-module is
    indecomposable and self-dual. Applying these results to
    $U=N^{\vee}$ and using duality (see \cite[Lemma 2.9.4 and
    Theorem 2.9.6]{RZ}), we obtain the decomposition.

    Using \eqref{it:lt3} together with $M_j^G =
    M_j^{(\bar\sigma-1)^{j-1}}$ for all $1\le j\le p$, we calculate
    \begin{align*}
        Z(T)\cap N &= \left(M_1\right)^{\Ic_i} \times
        \left(M_2^{(\bar\sigma-1)}\right)^{\Ic_2}\times \cdots
        \times \left(M_p^{(\bar\sigma-1)^{p-1}}\right)^{\Ic_p},\\
        Z(T)\cap T_{(i)} &=
        \left(M_i^{(\bar\sigma-1)^{i-1}}\right)^{\Ic_i} \times
        \cdots \times
        \left(M_p^{(\bar\sigma-1)^{p-1}}\right)^{\Ic_p}, \qquad 2\le
        i\le p,
    \end{align*}
    and $Z(T)\cap T_{(i)}=\{1\}$ for $i>p$.  We deduce that
    $t_i=\vert \Ic_i\vert$, $2\le i\le p$.

    For the case $i=1$, suppose first that $T$ is abelian of
    exponent $p$.  Then $Z(T)=Z(T)[p]=T$ and $T_{(2)}=\{1\}$.  By
    \eqref{it:lt3}, $N^{(\bar\sigma-1)}=\{1\}$, whence $\vert \Ic_i
    \vert=0$ for $i\ge 2$.  Therefore $t_1 = 1 + \vert \Ic_1 \vert$.
    Next suppose that $T$ is nonabelian. Then $Z(T)\subset N$.  We
    obtain $Z(T)[p]=Z(T)\cap N$ and so $t_1 = \dim_{\Fp} H^1 \left(
    (Z(T)\cap N)/(Z(T)\cap T_{(2)}), \Fp \right) = \vert\Ic_1\vert$.
    Finally, assume that $T$ is abelian and not of exponent $p$.
    Then $N=Z(T)[p]$ and $t_1 = \dim_{\Fp} H^1(N,\Fp) = \vert \Ic_1
    \vert$.

    (\ref{it:lt5}). For $\delta\in N$ we have $(\delta\sigma)^2 =
    [\sigma,\delta]\delta^2\sigma^2$ and, by induction,
    \begin{equation*}
        (\delta\sigma)^i =
        {[\underbrace{\sigma,[\sigma,\dots,[\sigma}_{i-1\ \text{times}},
        \delta]\cdots]]}^{\Tiny\begin{pmatrix}i\\ i\end{pmatrix}} \cdots
        [\sigma,[\sigma,\delta]]^{\Tiny\begin{pmatrix}
        i\\ 3\end{pmatrix}} [\sigma,\delta]^{\Tiny\begin{pmatrix}i\\
        2\end{pmatrix}} \delta^i\sigma^i.
    \end{equation*}
    Then $(N\sigma)^p= (\sigma^p)\cdot {}^{p-1}[\sigma, N]$, which
    by \eqref{it:lt2} and \eqref{it:lt3} may be written $\sigma^p
    \cdot T_{(p)}$. Replacing $\sigma$ with $\sigma^v$ for
    $(v,p)=1$, we conclude $T^p = \langle \sigma^p \rangle \cdot
    T_{(p)}$.

    (\ref{it:lt6}). Suppose that for some $v<u$ there is a
    commutative diagram \eqref{eq:cd} with a lift of $\bar\sigma$ in
    $H$ of order $p^2$.  Then $T\twoheadrightarrow H$ factors
    through $T/N^{(\bar\sigma-1)^{v}}$. But by \eqref{it:lt3},
    $T_{(v+1)} = N^{(\bar\sigma-1)^v}$ and by definition of $u$, we
    have $T^p\subset T_{(u)}\subset T_{(v+1)}$. Hence every lift of
    $\bar\sigma$ into $T/N^{(\bar\sigma-1)^v}$ is of order $p$, and
    the same holds for $H$.  We conclude that no commutative diagram
    as above with $\bar\sigma$ lifting to an element of order $p^2$
    exists for $v<u$.

    Now suppose that $u<p$ and consider $\sigma^p$.  By
    \eqref{it:lt5}, $T^p=\langle \sigma^p\rangle \cdot T_{(p)}$, and
    $T_{(p)}\subset T_{(u+1)}$.  We have $T^p\subset T_{(u)}$ and
    $T^p\not\subset T_{(u+1)}$ and so $\sigma^p\in T_{(u)}\setminus
    T_{(u+1)}$. Now $\sigma^p\in N$. By \eqref{it:lt3}, since
    $T_{(u)}=N^{(\sigma-1)^{u-1}}$ for $u\ge 2$, we deduce
    $\sigma^p\in N^{(\sigma-1)^{u-1}}\setminus N^{(\sigma-1)^u}$.
    From $[\sigma,\sigma^p]=1$ we obtain $\sigma^p\in N^G$.
    Therefore $\sigma^p \in (N^G\cap N^{(\sigma-1)^{u-1}}) \setminus
    (N^G \cap N^{(\sigma-1)^u})$.  We claim that there exists an
    $\Fp G$-submodule $M_u$ of $N$ such that $M_u^G=\langle
    \sigma^p\rangle$ and $N=M_u\times \tilde N$ for some $\Fp
    G$-submodule $\tilde N$ of $N$. Assume a factorization of $N$
    into cyclic $\Fp G$-submodules as in \eqref{it:lt4}.  Consider
    $w \in N$ such that $w^{(\sigma-1)^{u-1}} = \sigma^p$ and all
    components of $w$ lie in factors of dimension at least $u$. For
    at least one factor $M_u$, $\proj_{M_{u}} w$ generates $M_{u}$
    as an $\Fp G$-module. Write $N = M_u \times \tilde N$.  We
    obtain the commutative diagram
    \begin{equation*}
        \xymatrix{ 1 \ar[r]& N \ar[r] \ar@{>>}[d] & T \ar[r]
        \ar@{>>}[d] & G \ar[r] \ar[d]^{=} & 1 \\ 1 \ar[r] &
        M_u\ar[r] & T/\tilde N \ar[r] & G \ar[r] & 1}
    \end{equation*}
    in which a lift of $\bar\sigma$ is of order $p^2$.
\end{proof}

The following lemma follows easily from the definitions.
\begin{lemma}\label{le:t1}
    Let $T$ be a $T$-group with invariants $t_i$, $i=1, 2, \dots, p$,
    and $u$.
    \begin{enumerate}
        \item\label{it:tb1} $T$ is abelian if and only if $t_i=0$
        for all $2\le i\le p$.
        \item\label{it:tb2} $T$ is of exponent $p$ if and only if
        $u=p$ and $t_u=0$.
        \item\label{it:tb3} If $u<p$ then $t_u\ge 1$.
        \item\label{it:tb4} If $t_i=0$ for all $2\le i\le
        p$, then $t_1\ge 1$.
    \end{enumerate}
\end{lemma}

\begin{proof}[Proof of Proposition~\ref{pr:tc}]
    By Lemma~\ref{le:t1}(\ref{it:tb3},\ref{it:tb4}), \eqref{it:tc1}
    implies \eqref{it:tc2}.  Now suppose we are given cardinalities
    $t_i$, $i=1, 2, \dots, p$, and $u$ satisfying
    conditions~\eqref{it:tc2}.

    \emph{The case $u<p$}. Let $G$ be a group of order $p$ and
    $\Ic_i$, $i=1, 2, \cdots, p$ sets with cardinalities $\vert
    \Ic_i\vert$ satisfying $\vert \Ic_i\vert = t_i$ for $i\neq u$
    and $1+\vert \Ic_u\vert = t_u$. Set $N = X \times M_1^{\Ic_1}
    \times M_2^{\Ic_2} \times \cdots \times M_p^{\Ic_p}$, where
    $X\simeq M_u$ and $N$ is endowed with the product topology. Then
    $N$ is a pseudocompact $\Fp G$-module. (See
    \cite[page~443]{Br}.)  Define an action of $\Z_p$ on $N$ by
    letting a generator $\sigma$ of $\Z_p$ act via a generator of
    $G$, and form $N\rtimes \Z_p$ in the category of pro-$p$-groups.
    Now choose an $\Fp G$-module generator $x$ of $X$ and define
    $x_i=x^{(\sigma-1)^i}$ for $0\le i\le u$.  Since $(\sigma-1)$ is
    nilpotent of degree $u$ on $X$ we obtain $x_{u-1}\neq 1$ and
    $x_u=1$.  We set $R$ to be the closed subgroup $\langle
    \sigma^px_{u-1}\rangle \subset N\rtimes \Z_p$. Finally form $\G
    = (N\rtimes\Z_p) /R$ and set $\N$ to be the image of $N\rtimes
    \{1\}$ in $\G$. Since $\N\simeq N$ as pro-$p$-$G$ operator
    groups, we identify them.  By construction $\N$ is a maximal
    closed subgroup of $\G$ which is abelian of exponent $p$.
    Hence $\G$ is a $T$-group. Since the image of $\sigma$ in
    $\Gamma$ has order $p^2$, $\G$ is not of exponent $p$. From the
    decomposition of $N$, we obtain by
    Lemma~\ref{le:t}\eqref{it:lt4} that the invariants $t_i$ are as
    desired.  It remains only to show that $u$ is as given. By
    Lemma~\ref{le:t}\eqref{it:lt5} we have $\G^p = \langle
    x_{u-1}\rangle\cdot \G_{(p)}$.  From
    Lemma~\ref{le:t}(\ref{it:lt2},\ref{it:lt3}) we calculate that
    $x_{u-1}\in \G_{(u)}$ and $x_{u-1}\not\in \G_{(u+1)}$.  Hence
    $u$ is as desired.

    \emph{The case $u=p$} follows analogously.  Let $G$ be a group
    of order $p$ and $\Ic_i$ sets with cardinalities $\vert
    \Ic_i\vert$ satisfying $\vert \Ic_i\vert = t_i$, $2\le i\le p$;
    $\vert \Ic_1\vert = t_1$ if some $t_j\neq 0,\ 2\le
    j\le p$; and $1+\vert \Ic_1\vert = t_1$ if all $t_j=0$,
    $2\le j\le p$. Set $N = M_1^{\Ic_1}\times M_2^{\Ic_2} \times
    \cdots \times M_p^{\Ic_p}$, $\G=N\rtimes G$, $\N=N\rtimes
    \{1\}$, let $\sigma$ be a generator of $G$, and proceed as
    before.

    Now we show that $T$-groups are uniquely determined up to
    isomorphism by the invariants $t_i$ and $u$.  Let $T$ be an
    arbitrary $T$-group with invariants $t_i$, $i=1, 2, \dots, p$,
    and $u$, $N$ a maximal closed subgroup of $T$ that is abelian of
    exponent dividing $p$, and $G=T/N$.  From
    Lemma~\ref{le:t}(\ref{it:lt4}) the structure of $N$ as an $\Fp
    G$-module is determined up to isomorphism, and $T$ is an
    extension of $N$ by $G$.  Let $\sigma\in T\setminus N$.  We have
    $\sigma^p\in N^G$.  It remains only to determine the isomorphism
    class of $N$ as an $\Fp G$-module with a distinguished factor
    $X$ such that $\sigma^p\in X^G$.

    Suppose first that $u = p$, $t_1\ge 1$, and $t_i = 0$ for all $2
    \le i \le p$. From Lemma~\ref{le:t1}(\ref{it:tc1}), $T$ is
    abelian and so $T_{(2)}=\{1\}$.  Since $u=p$, $T$ has exponent
    $p$.  Then $T \simeq M_1^{\Ic_1} \times G$ and $t_1 = \vert
    \Ic_1 \vert + 1$.  Thus $T$ is determined by the invariants.
    Now suppose that $t_1 \ge 1$, $t_i = 0$ for all $2 \le i \le p$,
    and $u = 1$.  Again $T$ is abelian.  Since $u\neq p$, $T$ has
    exponent $p^2$.  Then $N = X \times \tilde N$, where $\sigma^p$
    generates an $\Fp G$-module $X$ isomorphic to $M_1$ and $\tilde
    N\simeq M_1^{\Ic'}$ with $\vert \Ic' \vert+1 = t_1$.

    Finally suppose $t_i \ne 0$ for some $i$ with $2\le i\le p$.
    Then $T$ is nonabelian, and by Lemma~\ref{le:t}(\ref{it:lt5}),
    $T^p= \langle\sigma^p\rangle \cdot T_{(p)}$.  From
    Lemma~\ref{le:t}\eqref{it:lt3} we obtain, for $u<p$, $\sigma^p
    \in ( N^G\cap N^{(\sigma-1)^{u-1}}) \setminus (N^G \cap
    N^{(\sigma-1)^u})$, while for $u=p$, we have $\sigma^p \in
    N^{(\sigma-1)^{p-1}}$.  If $u<p$ then by Proposition~\ref{pr:tc}
    we have $t_u\ge 1$, and using the same argument as in the proof
    of Lemma~\ref{le:t}\eqref{it:lt6}, $N$ contains a distinguished
    direct factor $X\simeq M_u$ such that $X^G=\langle
    \sigma^p\rangle$. We deduce that $N = X \times \tilde N$ where
    $\tilde N \simeq M_u^{\Ic'} \times \prod_{i\neq u} M_i^{\Ic_i}$
    for sets $\Ic_i$, $i\neq u$, and $\Ic'$ such that
    $\vert\Ic_i\vert = t_i$, $i\neq u$, and $1+\vert \Ic'\vert =
    t_u$. If $u=p$, we claim that without loss of generality we may
    assume that $\sigma^p=1$.  Since $\sigma^p\in
    N^{(\sigma-1)^{p-1}}$, let $\nu\in N$ such that
    $\sigma^p=\nu^{(\sigma-1)^{p-1}}$ and set $\tau =
    \sigma\nu^{-1}$.  Then $\tau\in T\setminus N$ and $\tau^p =
    \sigma^p (\nu^{-1})^{1+\sigma+\cdots+\sigma^{p-1}}=1$. Hence
    $T=N\rtimes G$.
\end{proof}

\begin{lemma}\label{le:freet}
    Suppose that $\G$ is a profinite group such that its maximal
    pro-$p$-quotient $\G(p)$ is a free pro-$p$-group of (positive
    and possibly infinite) rank $n$, and let $\N$ be a normal
    subgroup of $\G$ of index $p$.  Then the invariants of the
    $T$-group $\G/\N^p[\N,\N]$ are $t_1=1$, $t_i=0$ for $2\le i<p$,
    $t_p=n-1$ if $n<\infty$ and $t_p=n$ for $n$ an infinite
    cardinal, and $u=1$.
\end{lemma}

\begin{proof}
    Since $\G/\N^p[\N,\N]=\G(p)/\Phi(\N(p))$, we may assume without
    loss that $\G$ is a free pro-$p$-group. The result then follows
    from the analogue of the Kurosh subgroup theorem in the context
    of pro-$p$-groups.
\end{proof}

\begin{lemma}\label{le:freegf}
    Let $S$ be a free pro-$p$-group.  Then there exists a field $F$
    of characteristic $0$ such that $G_F\simeq S$.
\end{lemma}

\begin{proof}
    First let $F_0$ be any algebraically closed field of
    characteristic $0$ with cardinality greater than or equal to $d
    = \dim_{\Fp} H^1(S,\Fp)$.  Set $F_1:=F_0(t)$.  By \cite{Do},
    $G_{F_1}$ is a free profinite group, and let $P$ denote a
    $p$-Sylow subgroup of $G_{F_1}$.  By \cite[Corollary~7.7.6]{RZ},
    $P$ is a free pro-$p$-group.  Let $F_2$ be the fixed field of
    $P$.  The classes in $F_1^\times/F_1^{\times p}$ of the set $A$
    of linear polynomials $t-c$, $c\in F_0$, are linearly
    independent over $\Fp$.  Choose a subset of $A$ of cardinality
    $d$, and let $V$ be the vector subspace of $F_1^\times/
    F_1^{\times p}$ generated by the classes of the elements of $A$.
    Since $([F_2:F_1],p)=1$, $V$ injects into
    $F_2^{\times}/F_2^{\times p}$.  Let $W$ denote this image.  Now
    let $F$ be a maximal algebraic field extension of $F_2$ such
    that $W$ injects into $F^\times/F^{\times p}$.  By maximality,
    the image $i(W)$ of $W$ in $F^{\times}/F^{\times p}$ is in fact
    $F^{\times}/F^{\times p}$.  The rank of $G_F$ is then
    $\dim_{\Fp} H^1(G_F,\Fp) = \dim_{\Fp} V = d$.
\end{proof}

\begin{proof}[Proof of Theorem~\ref{th:t}]
    \emph{The case $p=2$}. Let $u\in \{1,2\}$, $t_1$, and $t_2$ be
    invariants of a $T$-group $T$.  By Proposition~\ref{pr:tc},
    $t_1\ge 1$ if $u=1$, and if $u=2$ then either $t_1\ge 1$ or
    $t_2\ge 1$.

    \emph{Case 1}: $T$ is not of exponent $2$.  By
    Lemma~\ref{le:t1}(\ref{it:tb2}), either $u=1$ or $t_2\ge 1$. Let
    $G$ be a group of order $2$ and $N=M_1^{\Ic_1}\times
    M_2^{\Ic_2}$ for sets $\Ic_1$ and $\Ic_2$ satisfying
    $\vert\Ic_1\vert=t_1$ and $\vert\Ic_2\vert=t_2$, and let
    $M=N^\vee$.  From \cite[Corollary~2]{MSw2}, there exists $E/F$
    with $\chr F\neq 2$ and $\Gal(E/F)\simeq G$ such that
    $H^1(G_E,\Ft)\simeq E^\times/E^{\times 2} \simeq M$ as $\Ft
    G$-modules if and only if there exist $\Upsilon\in \{0,1\}$ and
    cardinalities $d$ and $e$ such that $t_1+1=2\Upsilon+d$;
    $t_2+\Upsilon=e$; if $\Upsilon=0$ then $d\ge 1$; and if
    $\Upsilon=1$ then $e\ge 1$. Moreover, $-1\in N_{E/F}(E^\times)$
    if and only if $\Upsilon=1$.  Finally, by \cite[proof of
    Theorem~1]{MSw2}, we may choose $E/F$ such that $G_F$ is a
    pro-$2$-group.

    If $u=1$ then set $\Upsilon=1$ and $e=t_2+\Upsilon$.  Since
    $t_1\ge 1$ we may choose $d\ge 0$ such that $2\Upsilon+d=t_1+1$.
    Then $e\ge 1$ and the conditions for $E/F$ with
    $G_E/G_E^{(2)}\simeq M^\vee \simeq N$ are satisfied. Since
    $\Upsilon=1$, $-1\in N_{E/F}(E^\times)$, and by Albert's
    criterion \cite{A}, $E/F$ embeds in a cyclic extension $E'/F$ of
    degree $4$.  Let $\Gal(E/F)=\langle \bar\sigma\rangle$. We have
    the commutative diagram
    \begin{equation*}
        \xymatrix{ 1 \ar[r]& G_E/G_E^{(2)} \ar[r] \ar@{>>}[d] &
        G_F/G_E^{(2)} \ar[r] \ar@{>>}[d] & G \ar[r] \ar[d]^{=} & 1
        \\ 1 \ar[r] & M_1 \ar[r] & H \ar[r] & G \ar[r] & 1}
    \end{equation*}
    in which $\bar\sigma$ lifts to an element of order $4$.  By
    Lemma~\ref{le:t}\eqref{it:lt6}, $u=1$ for $T_{E/F}$. By
    Lemma~\ref{le:t1}\eqref{it:tb2}, $T_{E/F}$ is not of exponent
    $2$.  Using Lemma~\ref{le:t}\eqref{it:lt4}, because
    $G_E/G_E^{(2)}\simeq N\simeq M_1^{\Ic_1} \times M_2^{\Ic_2}$,
    the invariants $t_1$ and $t_2$ of $T_{E/F}$ match those of $T$.
    By Proposition~\ref{pr:tc}, $T\simeq T_{E/F}$.

    If $u=2$ then $t_2\ge 1$.  We take $\Upsilon=0$, $d=t_1+1\ge 1$,
    and $e=t_2$ and obtain an extension $E/F$ as before. Since
    $\Upsilon=0$, $-1\not\in N_{E/F}(E^\times)$ and so by \cite{A},
    $E/F$ does not embed in a cyclic extension $E'/F$ of degree $4$.
    Let $\Gal(E/F)=\langle \bar\sigma\rangle$. There is no
    commutative diagram
    \begin{equation*}
        \xymatrix{ 1 \ar[r]& G_E/G_E^{(2)} \ar[r] \ar@{>>}[d] &
        G_F/G_E^{(2)} \ar[r] \ar@{>>}[d] & G \ar[r] \ar[d]^{=} & 1
        \\ 1 \ar[r] & M_1 \ar[r] & H \ar[r] & G \ar[r] & 1}
    \end{equation*}
    in which $\bar\sigma$ lifts to an element of order $4$.  By
    Lemma~\ref{le:t}\eqref{it:lt6}, $u=2$ for $T_{E/F}$. Because
    $t_2\ge 1$, $G_E/G_{E}^{(2)}$ contains an $\Ft G$-submodule
    isomorphic to $M_2$ whence $T_{E/F}$ is nonabelian.  By
    Lemma~\ref{le:t}\eqref{it:lt4} and the isomorphism
    $G_E/G_E^{(2)}\simeq N\simeq M_1^{\Ic_1} \times M_2^{\Ic_2}$ we
    deduce that the invariants $t_1$ and $t_2$ of $T_{E/F}$ match
    those of $T$.  By Proposition~\ref{pr:tc}, $T\simeq T_{E/F}$.

    \emph{Case 2}: $T$ has exponent $2$. By
    Lemma~\ref{le:t1}(\ref{it:tb2}), $u=2$ and $t_2=0$, and so
    $t_1\ge 1$.  Let $N=M_1^{\Ic_1}$ for $\Ic_1$ satisfying
    $\vert\Ic_1\vert+1=t_1$. Set $\Upsilon=0$, $d=t_1+1$, and
    $e=t_2=0$.  Then $d\ge 1$ and there exists $E/F$ with $\chr
    F\neq 2$ such that $H^1(G_E,\Ft)\simeq M$ and $-1\not\in
    N_{E/F}(E^\times)$.  As before, $u=2$ for $T_{E/F}$, and by
    Lemma~\ref{le:t}\eqref{it:lt4}, $t_2=0$ for $T_{E/F}$. By
    Lemma~\ref{le:t1}(\ref{it:tb2}), $T_{E/F}$ has exponent $2$, and
    by Lemma~\ref{le:t}\eqref{it:lt4}, $t_1$ is the correct
    invariant for $T_{E/F}$. By Proposition~\ref{pr:tc}, $T\simeq
    T_{E/F}$.

    \emph{The case $p>2$}.  First we characterize those $T$-groups
    occurring as $T_{E/F}$ for fields $F$ such that the maximal
    pro-$p$-quotient $G_F(p)$ of the absolute Galois group $G_F$ is
    free pro-$p$.  Lemma~\ref{le:freet} tells us that for such $F$
    and an $E/F$ of degree $p$, the invariants of $T_{E/F}$ are
    $t_1=1$, $t_i=0$ for $2\le i<p$, and $u=1$, and that the rank of
    $G_F(p)$ is one more than the invariant $t_p$.  Now suppose that
    $T$ is a $T$-group with invariants $t_1=1$, $t_i=0$ for $2\le
    i<p$, and $u=1$.  By Lemma~\ref{le:freegf} there exists
    $F$ such that $G_F$ is a free pro-$p$-group of rank $t_p+1$.
    Letting $\N$ be any maximal closed subgroup of $G_F$ and
    $E=\Fix(\N)$, Lemma~\ref{le:freet} and Proposition~\ref{pr:tc}
    give $T\simeq T_{E/F}$. Therefore the $T$-groups which occur as
    $T_{E/F}$ for fields $F$ with free maximal pro-$p$-quotient
    $G_F(p)$ are precisely those for which $t_1=1$, $t_i=0$ for
    $2\le i<p$, and $u=1$.

    Now we characterize which of the remaining $T$-groups occur as
    $T_{E/F}$ for cyclic field extensions $E/F$ of degree $p$ for
    $F$ a field such that either $\chr F=p$ or $\xi_p\in F$.  By
    Witt's Theorem, $\chr F\neq p$. Hence we consider only fields
    $F$ with $\chr F\neq p$ and $\xi_p\in F$.  For a cyclic
    extension $E/F$ of degree $p$, consider the
    $\Fp\Gal(E/F)$-module $M_{E/F}:=H^1(G_E,\Fp)\simeq
    E^\times/E^{\times p}$, and let $G$ be an abstract group of
    order $p$.  Since the particular isomorphism $G\simeq \Gal(E/F)$
    does not alter the isomorphism class of $M_{E/F}$ as an $\Fp
    G$-module, we may consider all such modules $\Fp G$-modules.
    (See \cite{MSw2}.) By \cite[Corollary~1]{MSw2}, $M\simeq
    M_{E/F}$ as $\Fp G$-modules for suitable $E/F$ with $G\simeq
    \Gal(E/F)$ if and only if $M = M_1^{\Ic_1} \oplus M_2^{\Ic_2}
    \oplus M_p^{\Ic_p}$, where the cardinalities $j_1=\vert
    \Ic_1\vert$, $j_2=\vert \Ic_2\vert$, and $j_p=\vert\Ic_p\vert$
    satisfy the following conditions: $j_1+1= 2\Upsilon+d$,
    $j_2=1-\Upsilon$, and $j_p+1=e$ for some cardinalities $d$, $e$
    and $\Upsilon\in \{0,1\}$ where $d\ge 1$ if $\Upsilon=0$ and
    $e\ge 1$.  Moreover, $\Upsilon=1$ if and only if $\xi_p\in
    N_{E/F}(E^\times)$. Finally, by \cite[proof of Theorem~1]{MSw2},
    we may choose $E/F$ such that $G_F$ is a pro-$p$-group.

    We observe that the constraints on $j_1$, $j_2$, and $j_p$ are
    then $j_1\ge \Upsilon$ and $j_2=1-\Upsilon$.  Now by duality,
    $H^1(G_E,\Fp)^\vee \simeq G_E/G_E^{(2)}$, and since cyclic $\Fp
    G$-modules are self-dual, we may derive conditions on the
    topological $\Fp G$-module $G_E/G_E^{(2)}$ occurring as maximal
    closed subgroups of $T$-groups $G_F/G_E^{(2)}$, as follows. Set
    $G=G_F/G_E$.

    First we relate $\Upsilon$ and the invariant $u$. We claim that
    for any $T$-group $G_F/G_E^{(2)}$, we have $u\le 2$.  Write
    $E=F(\root{p}\of{a})$ for some $a\in F^\times$.  Let $[e]\in
    E^\times/E^{\times p}$ denote the class of $e\in E^\times$. Then
    $X=\langle [\root{p}\of{a}], [\xi_p]\rangle$ is a cyclic $\Fp
    G$-submodule of $M$ and is isomorphic to $M_i$ for some $i \in
    \{ 1,2 \}$. By equivariant Kummer theory (see \cite{W}),
    $L=E(\root{p^2} \of{a},\xi_{p^2})$ is a Galois extension of $F$.
    Moreover, $G_F/G_L$ is a homomorphic image of $G_F/G_E^{(2)}$,
    since $L/E$ is an elementary abelian extension. Then
    $L/F(\xi_{p^2})$ is Galois with group $\Z/p^2\Z$, and for any
    $\sigma\in G_F\setminus G_E$, the restriction $\sigma_L\in
    \Gal(L/F)$ restricts to a generator. Hence $1\neq \sigma^p\in
    \Gal(L/E)$, and therefore $1\neq \sigma^p\in G_E/G_E^{(2)}$.  We
    have a commutative diagram
    \begin{equation}\label{eq:cd2}
        \vcenter{\xymatrix{ 1 \ar[r]& G_E/G_E^{(2)} \ar[r]
        \ar@{>>}[d] & G_F/G_E^{(2)} \ar[r] \ar@{>>}[d] & G \ar[r]
        \ar[d]^{=} & 1 \\ 1 \ar[r] & M_i \ar[r] & G_F/G_L \ar[r] & G
        \ar[r] & 1}}
    \end{equation}
    and so by Lemma~\ref{le:t}(\ref{it:lt6}) we deduce that $u\le
    i\le 2$.

    Now we claim that $\Upsilon=1$ if and only if $u=1$. We have
    that $\Upsilon=1$ if and only if $\xi_p\in N_{E/F}(E^\times)$.
    By \cite{A}, $\Upsilon=1$ if and only if $E/F$ embeds in a
    cyclic extension of $F$ of degree $p^2$, if and only if there
    exists a closed normal subgroup $\tilde N\subset G_E$ such that
    $G_F/\tilde N\simeq \Z/p^2\Z$.  Any such closed normal subgroup
    must contain $G_E^{(2)}$.  Hence we deduce that $\Upsilon=1$ if
    and only if there exists a commutative diagram \eqref{eq:cd2}
    with $i=1$ in which nontrivial elements of $G$ lift to elements
    of order $p^2$.  By Lemma~\ref{le:t}\eqref{it:lt6}, $\Upsilon=1$
    if and only if $u=1$.

    We have therefore shown that $u\le 2$ and $u=2-\Upsilon$.
    Translating the remaining conditions on $j_1$ and $j_2$, we see
    that $j_1\ge 2-u$ and $j_2=u-1$.  Now by
    Lemma~\ref{le:t}\eqref{it:lt4}, $t_2=j_2$, $t_i=0$ for $3\le
    i\le p-1$, and $t_p=j_p$.  Moreover, $t_1=j_1$ if $T$ is not
    abelian of exponent $p$.  By Lemma~\ref{le:t1}\eqref{it:tb2},
    $T$ is of exponent $p$ if and only if $u=p$.  But we have shown
    that $u\le 2<p$, whence $T$ is not of exponent $p$ and we have
    $t_1=j_1$. By Proposition~\ref{pr:tc}, if $u=1$ then $t_1\ge 1$,
    and since $t_2=u-1$ the conditions for applying
    \cite[Corollary~1]{MSw2} are valid. Hence a $T$-group $T_{E/F}$
    with prescribed invariants subject to conditions (1) exists.
\end{proof}

\begin{proof}[Proof of Theorem~\ref{th:1}]
    Suppose $\G$ is a pro-$p$-group with maximal closed subgroup
    $\N$, and let $T=\G/\N^{(2)}$, $N=\N/\N^{(2)}$ and $G =
    \Gamma/\Delta$.  Write $\bar\sigma$ and $\bar\tau$ for the
    images of $\sigma$ and $\tau$, respectively.

    (1).  Since $\N$ surjects onto $N$ we have that $\bar\sigma
    \not\in N$, $\bar\tau\in N$, ${}^e[\bar\sigma, \bar \tau]
    \not\in {}^{p-1}[\bar\sigma, N]$, and ${}^{e+1}[\bar\sigma,
    \bar\tau] =1$. By Lemma~\ref{le:t}(\ref{it:tb2},\ref{it:tb3}) we
    have $1\neq {}^e[\bar\sigma,\bar\tau] \in T_{(e+1)} \setminus
    T_{(p-1)}$. We deduce from ${}^{e+1}[\bar\sigma,\bar\tau]=1$
    that ${}^e[\bar\sigma, \bar\tau]\in Z(T)$.  We obtain that some
    invariant $t_{i}\neq 0$ for $3\le i<p$. Now if $\G=G_F$ for some
    $F$, then setting $E=\Fix(\N)$ we obtain $T=T_{E/F}$,
    contradicting Theorem~\ref{th:t}.  Now assume additionally that
    there exists a closed normal subgroup $\Lambda\subset \N$ of
    $\G$ such that $\G/\Lambda\simeq \Z/p^2\Z$ and $e=1$.  Let
    $\tilde\sigma$ denote the image of $\sigma$ in $T$.  We have a
    commutative diagram \eqref{eq:cd} with $v=1$ and an image of
    $\tilde\sigma$ in $H$ of order $p^2$. By
    Lemma~\ref{le:t}\eqref{it:lt6}, $u=1$ for $T$ is equal to $1$.
    As before some $t_i\neq 0$, $2\le i<p$. Again by
    Theorem~\ref{th:t} we are done.

    (2). We proceed as before, obtaining two elements
    $[\bar\sigma,\bar\tau_i]$, $i=1, 2$, which generate distinct
    subgroups of $Z(T)\cap T_{(2)}$ with trivial intersection with
    $T_{(p)}$.  We deduce that the sum of $t_i$, $2\le i<p$, is at
    least two, and we apply Theorem~\ref{th:t}.

    (3). We obtain $\bar\sigma\not\in N$, $\bar\tau\in N$, and
    $\bar\sigma^p \in {}^{2}[\bar\sigma, N]$. By
    Lemma~\ref{le:t}(\ref{it:tb2},\ref{it:tb3}) we have
    $\bar\sigma^p\subset T_{(3)}$, and from
    Lemma~\ref{le:t}(\ref{it:lt5}) we deduce that $T^p\subset
    T_{(3)}$.  Hence $u\ge 3$, and we apply Theorem~\ref{th:t}.
\end{proof}

\begin{proof}[Proof of Corollary to Theorem~\ref{th:1}]
    Let $\N$ be the closed subgroup of $\G$ generated as a normal
    subgroup by $\sigma_1^p$, $\sigma_2$, and $\sigma_i$ for $i\in
    \Ic$, and let $\Lambda$ be the closed subgroup of $\G$ generated
    as a normal subgroup by $\sigma_1^{p^2}$, $\sigma_2$, and
    $\sigma_i$ for $i\in \Ic$.  Set also $T = \Gamma/\Delta^{(2)}$.
    Examining the quotient of $\G$ obtained by trivializing
    $\sigma_2$ and each $\sigma_i$, $i\in \Ic$, we see that $\N$ is
    maximal and $\G/ \Lambda \simeq \Z/p^2\Z$. Now let $e=f-1$.  By
    passing from $\G$ to $T$ using bars for denoting images of
    elements of $\G$ in $T$, we see that
    ${}^e[\bar\sigma_1,\bar\sigma_2] \notin
    {}^{p-1}[\bar\sigma_1,\bar\Delta]$ in $T$. On the other hand,
    ${}^{e+1}[\sigma_1,\sigma_2]={}^f[\sigma_1,\sigma_2]=1$ in $T$.
    By Theorem~\ref{th:1}(1), $\G$ is not an absolute Galois group.
\end{proof}

It is a natural question to ask whether the maximal closed subgroup
$N$ in the definition of $T$-group is unique.  It is not difficult
to show that if $T$ is a $T$-group with a maximal closed subgroup
$N$ which is abelian of exponent dividing $p$, then unless $T$ is
either itself abelian of exponent $p$ of order greater than $p$, or
isomorphic to the direct product of a Heisenberg group of order
$p^3$ and an abelian group of exponent dividing $p$, then $N$ is a
characteristic subgroup of $T$. Since these exceptional $T$-groups
have invariant $u=p$, we have by Theorem~\ref{th:t} that for $p>2$,
each $T_{E/F}$ has a unique index $p$ elementary abelian subgroup.

\section{Analogues of Theorems of Artin-Schreier and
Becker}\label{se:asb}

Recall that by Artin-Schreier there are no elements of order $p$ in
the absolute Galois groups of a field $F$, unless $p=2$ and $F$ is
formally real. Becker proved that the same holds for maximal
pro-$p$-quotients of absolute Galois groups \cite{Be}. In this
section we show that in certain small quotients of absolute Galois
groups, there are no non-central elements of order $p$ unless $p=2$
and the base field $F$ is formally real.

Let $p$ be a prime and $F$ a field with $\xi_p\in F$.  We define the
following fields associated to $F$:
\begin{itemize}
    \item $F^{(2)}$: the compositum of all cyclic extensions of $F$
    of degree $p$
    \item $F^{(3)}$: the compositum of all cyclic extensions of
    $F^{(2)}$ of degree $p$ which are Galois over $F$
    \item ${F^{(2)}}^{(2)}$: the compositum of all cyclic extensions
    of $F^{(2)}$ of degree $p$.
\end{itemize}
Then we set $W_F = \Gal(F^{(3)}/F)$ and $V_F = \Gal({F^{(2)}}^{(2)}
/F)$.  Observe that for a cyclic extension $E/F$ of degree $p$,
$T_{E/F} =\Gal(E^{(2)}/F)$.

\begin{theorem}\label{th:2}
    Let $p$ be prime and $F$ a field with $\xi_p\in F$.  The
    following are equivalent:
    \begin{enumerate}
        \item\label{it:w2} There exists $\sigma\in V_F\setminus
        \Phi(V_F)$ of order $p$
        \item\label{it:w1} There exists $\sigma\in W_F\setminus
        \Phi(W_F)$ of order $p$
        \item\label{it:w0a} There exists $\sigma\in V_F \setminus
        \Phi(V_F)$ such that for each cyclic $E/F$ of degree $p$
        its image $\sigma_{E/F} \in T_{E/F}$ has order at most $p$.
        \item\label{it:w3} $p=2$ and $F$ is formally real.
    \end{enumerate}
    If these conditions hold, the elements whose square roots are
    fixed by $\sigma$ form an ordering of $F$.
\end{theorem}

\begin{proof}[Proof of Theorem~\ref{th:2}]
    $(1)\implies (2)$ and $(1)\implies (3)$. Observe first that
    $F^{(3)}\subset {F^{(2)}}^{(2)}$.  Hence we have a natural
    surjection $V_F \twoheadrightarrow W_F$. Assume that $\sigma \in
    V_F \setminus \Phi(V_F)$ has order $p$. Then its image in $W_F
    \setminus \Phi(W_F)$ also has order $p$. Moreover, for any
    cyclic extension $E/F$ of degree $p$, the image $\sigma_{E/F}$
    of $\sigma$ in $T_{E/F}$ has order at most $p$, as follows.
    Since $E^{(2)}$ is contained in ${F^{(2)}}^{(2)}$ we see that
    $\sigma_{E/F}^p= 1$ in $T_{E/F}$.

    $(2)\implies (4)$.  Suppose that $\sigma\in W_F\setminus
    \Phi(W_F)$ has order $p$.  As in the proof of
    Theorem~\ref{th:t}, let $a\in F^\times\setminus F^{\times p}$ be
    arbitrary such that $F(\root{p}\of{a})\not\subset F(\xi_{p^2})$.
    Set $K_a:=F(\root{p}\of{a},\xi_{p^2})\subset F^{(2)}$.  Then
    $[K_a:F(\xi_{p^2})]=p$ and $L_a := F (\root{p^2}\of{a},
    \xi_{p^2})$ is a cyclic extension of degree $p^2$ of
    $F(\xi_{p^2})$.  Moreover, $L_a/F$ is a Galois extension of $F$
    and $L_a \subset F^{(3)}$ since $\Gal(L_a/F)$ is a central
    extension of degree $\le p$ of $\Gal(K_a/F)$. Now if
    $\sigma(\root{p}\of{a}) \ne \root{p}\of{a}$ then $\sigma$ has
    order $p^2$ in $W_F$. Hence $\sigma$ fixes all $\root{p}\of{a}$
    for $a\in F^\times$ with $F\neq F(\root{p}\of{a})\not\subset
    F(\xi_{p^2})$.  Since $\sigma\not\in \Phi(W_F)$, there must
    exist a cyclic extension $E\subset F^{(2)}$ of $F$ which is not
    fixed by $\sigma$.  We deduce that $\xi_{p^2}\not\in F$ and so
    $E=F(\xi_{p^2})$.  If $p$ is odd then $L:=F(\xi_{p^3})$ is a
    cyclic extension of $F$ of degree $p^2$ and $\sigma$ restricts
    to a generator of $\Gal(L/F)\simeq \Z/p^2\Z$, again a
    contradiction, whence $p=2$.  Now let $s \in W_F \setminus
    \Phi(W_F)$ be an element of order $2$. By \cite[proof of
    Theorem~2.7]{MSp1}, the elements of $F$ whose square roots are
    fixed by $s$ form an ordering of $F$.  Hence $(2)\implies (4)$.

    $(4)\implies (1)$.  Suppose that $(4)$ holds.  By
    \cite[Satz~3]{Be}, there exists an ordering of $F$ whose square
    roots are fixed by some element $s$ of order $2$ in $G_F(2)$.
    Then the restriction of $s$ to ${F^{(2)}}^{(2)}$ is the required
    element $\sigma \in V_F \setminus \Phi(V_F)$ of order $2$.
    Hence $(4) \implies (1)$.

    $(3)\implies (2)$.  Let $\sigma \in V_F \setminus \Phi(V_F)$
    such that for each cyclic extension $E/F$ of degree $p$ the
    image $\sigma_{E/F}$ of $\sigma$ in $T_{E/F}$ has order at most $p$. Let
    $E = F(\root{p} \of{a})$ such that $\sigma$ acts nontrivially on
    $\root{p}\of {a}$ and $L_a = F(\root{p^2}\of{a},\xi_{p^2})$.  As
    above, since the restriction of $\sigma_{E/F}$ to $L_a \subset
    E^{(2)}$ is not of order $p^2$, we deduce that $p=2$ and
    $\sqrt{-1}\not\in F$. Hence $F$ is not quadratically closed.  If
    $F$ is real closed then there exists precisely one extension
    $E/F$ of degree $2$, namely $F^{(2)}$, and $F^{(2)}=F^{(3)}=
    {F^{(2)}}^{(2)}$, whence $W_F=T_{E/F}=V_F$. Otherwise, $F^{(3)}$
    is a compositum of extensions $L/F$ such that $\Gal(L/F)$ is
    either a cyclic group of order $4$ or a dihedral group of order
    $8$. (See \cite[Corollary~2.18]{MSp2}.) On the other hand, each
    such $L$ lies in $E^{(2)}$ for a suitable quadratic extension
    $E/F$: each $L$ may be obtained as a Galois closure of
    $E(\sqrt{\gamma})$ for some $[E:F] = 2$ and $\gamma \in
    E^\times$.  Therefore the restrictions $\sigma_{L/F}$ of
    $\sigma$ to the extensions $L/F$ have order $\le 2$, and so the
    restriction of $\sigma_{F^{(3)}/F}$ of $\sigma$ to $F^{(3)}/F$
    has order $2$. Hence $(3)\implies (2)$.
\end{proof}

\section{Pro-$p$-Groups That are Not Absolute Galois
Groups}\label{se:fam}

\begin{theorem}\label{th:fam}
    Let $p>3$ be prime.  There exists a pro-$p$-$\Z_p$ operator
    group $\Omega$ such that no group of the form $\G := ((\Omega
    \star \Sigma)\rtimes \Z_p)/\Ec$, where $\Sigma$ is any
    pro-$p$-group with trivial $\Z_p$-action, and $\Ec$ is any
    normal closed subgroup of $(\Omega\star\Sigma)\rtimes \Z_p$ such
    that $\Ec \subset ((\Omega\star\Sigma)\rtimes p\Z_p)^{(3)}$, is
    an absolute Galois group.
\end{theorem}
\noindent Here $\star$ denotes the free product in the category of
pro-$p$-groups.  (Recall that $R^{(n)}$ denotes the $n$th term of
the $p$-descending series of a pro-$p$-group $R$; see
section~\ref{se:t}.)

The $\Omega$ of the theorem is that of the following proposition.
Recall that for a pro-$p$-group $\G$, the \emph{decomposable
subgroup} of $H^2(\G,\Fp)$ is defined to be the subgroup generated
by cup products of elements of $H^1(\G,\Fp)$: $H^2(\G,\Fp)^\dec =
H^1(\G,\Fp).H^1(\G,\Fp)$. We say that $H^2(\G,\Fp)$ is
\emph{decomposable} if $H^2(\G,\Fp)= H^2(\G,\Fp)^\dec$.

\begin{proposition}\label{pr:const}
    Let $p>3$ be prime and $C$ be a cyclic group of order $p$. There
    exists a torsion-free pro-$p$-$C$ operator group $\Omega$ such
    that as $\Fp C$-modules, $H^1(\Omega,\Fp) \simeq M_p$ and
    $H^2(\Omega,\Fp) = H^2(\Omega,\F_p)^\dec \simeq M_{p-1}\oplus
    \frac{(p-3)}{2} M_p$.
\end{proposition}
\noindent Here we use $M_i$ to denote the unique cyclic $\Fp
C$-module of dimension $i$, and we write the action of $\Fp C$
multiplicatively.

\begin{proof}
    Let $D=\langle g \mid g^p=1 \rangle$ be a cyclic group of order
    $p$, and let $\Z_p D$ be the $p$-adic group ring, written
    multiplicatively as $\GG$, where the element $\bar g_i$ of $\GG$
    corresponds to the element $g^i$ of $\Z_p D$. We interpret the
    suffixes mod $p$. Now let $C=\langle \sigma \mid
    \sigma^p=1\rangle$ be another group of order $p$, acting on
    $\GG$ via $\sigma(\bar g_i)=\bar g_{i-1}$. In this way $\GG$ and
    $H^1(\GG,\Z_p)$ are free $\Z_p C$-modules, and $H^1(\GG,\Z_p)$
    has a topological generating set $y_0, y_1, \dots, y_{p-1}$ dual
    to $\bar g_0, \bar g_1, \dots, \bar g_{p-1}$.  Observe that
    $\sigma(y_i)=y_{i+1}$. Next let $\H=\Z_p$ be a trivial $\Z_p
    C$-module with generator $h$, and let $z\in H^1(\H,\Z_p)$ be
    dual to $h$.  We define a nonsplit central extension $\Omega$ of
    $\H$ by $\GG$ as follows. The group $H^2(\GG,\H) = \bigwedge^2
    H^1(\GG,\H)$ is a free $\Z_p C$-module of rank $p(p-1)/2$ with
    free generators $y_0y_j$ for $1\le j\le (p-1)/2$. Consider the
    element
    \begin{equation*}
        y := (1+\sigma+\dots+\sigma^{p-1})y_0y_1 = y_0y_1+y_1y_2+\dots
        +y_{p-1}y_0.
    \end{equation*}
    Let $\Omega$ be the central extension of $\H$ by $\GG$
    corresponding to $y\in H^2(\GG,\H)$.

    The group $\Omega$ is a torsion-free nilpotent pro-$p$-group of
    Hirsch length $p+1$. The standard correspondence of group
    extensions with $H^2$ gives us that for suitable representatives
    $g_i\in \Omega$ with images $\bar g_i\in \GG$, we have the
    relations $[g_i,g_j] = h$ for $j=i+1$; $[g_i,g_j] = 1$ for
    $j\neq i+1$ and $i\neq j+1$; and $[g_i,h] = 1$ for all $i$. The
    action of $C$ on $\GG$ may be extended to $\Omega$ by
    $\sigma(g_i)=g_{i-1}$, $0\le i< p$, and $\sigma(h)=h$.

    The $E_2$ page of the spectral sequence $H^s(\GG,H^t(\H,\Fp))
    \Rightarrow H^{s+t}(\Omega,\Fp)$ is generated by anticommuting
    elements $\bar y_i\in E_2^{1,0}$ and $\bar z\in E_2^{0,1}$, and we
    have $d_2^{0,1}(\bar z)=\bar y$, where $\bar y$ is the reduction
    mod $p$ of $y$. Observe that $d_2^{0,1}$ is injective on
    $E^{0,1}=H^1(\H,\F_p)$. By the five-term exact sequence,
    $H^1(\Omega,\Fp)\simeq H^1(\GG,\Fp)$ as $\Fp C$-modules and so
    $H^1(\Omega,\Fp)\simeq M_p$.

    Now consider $H^2(\Omega,\Fp)$. We claim first that $d_2^{1,1}$
    is injective. It is enough to show that $d_2^{1,1}(\bar z.\bar
    x)=0$ implies that $\bar x=0$ for $\bar x\in H^1(\GG,\Fp)$.
    Write $\bar x=\sum_{i=0}^{p-1} c_i \bar y_i$, $c_i\in \Fp$.
    Since $p>3$, the set of elements $\bar y_i.\bar y_{i+2}.\bar
    y_{i+3}$, $0\le i<p$, is $\Fp$-independent and may be expanded
    to an $\Fp$-basis of $E_2^{3,0}=H^3(\GG,\Fp)$ consisting of
    threefold products of $\bar y_i$. Consider the coefficient
    $\beta_i$ of $\bar y_i.\bar y_{i+2}.\bar y_{i+3}$ in an
    expansion of $d_2^{1,1}(\bar z.\bar x)$. Since $p>3$, the only
    consecutive pair of indices in $\bar y_i.\bar y_{i+2}.\bar
    y_{i+3}$ is $\{i+2,i+3\}$. Hence $\beta_i$ depends only on
    $c_i$. If all $\beta_i=0$, then each $c_i=0$ and therefore
    $d_2^{1,1}$ is injective. Now since $d_2^{1,1}$ is injective,
    $E_3^{1,1} = E_\infty^{1,1}=0$, and because $E_2^{0,2}=0$ we have
    $E_\infty^{0,2}=0$. Next observe that $\bar y =
    (\sigma-1)^{p-1}\bar y_0.\bar y_1$ and so $\bar y\in
    H^2(\GG,\Fp)^C$.  Since $E^{2,0}=H^2(\GG,\Fp)$ is the direct sum
    of free $\Fp C$-modules $M(0,i)\simeq M_p$ on generators $\bar
    y_0.\bar y_i$, $1\le i\le (p-1)/2$, we deduce that
    $E_\infty^{2,0}=E_3^{2,0}=(M(0,1)/M(0,1)^C) \oplus
    \frac{(p-3)}{2} M(0,i)$. Thus $H^2(\Omega,\Fp) \simeq M_{p-1}
    \oplus \frac{(p-3)}{2} M_p$.  Finally, since
    $H^1(\Omega,\Fp)\simeq H^1(\GG,\Fp)$ and $H^2(\Omega,\Fp)$ is a
    quotient of $E_2^{2,0}=H^2(\GG,\Fp)= \bigwedge^2 H^1(\GG,\Fp)$,
    $H^2(\Omega,\Fp)$ is decomposable as well.
\end{proof}

Given a free pro-$p$-group $V$ and a pro-$p$-group $\N$, we say that
a surjection $V\twoheadrightarrow \N$ is a \emph{minimal
presentation} of $\N$ if $\inf:H^1(\N,\Fp)\to H^1(V,\Fp)$ is an
isomorphism.

\begin{proposition}\label{pr:a2}
    Let $\Delta$ be a pro-$p$-group, $V$ a free pro-$p$-group, and
    $1 \to R \to V \to \N \to 1$ a minimal presentation of $\N$.
    Then we have natural maps
    \begin{enumerate}
        \item\label{it:n1} $H^2(\N,\Fp)\simeq (R/(R^p[R,V]))^\vee$
        \item\label{it:n2} $H^2(\N,\Fp)^\dec \simeq (R/(R\cap
        V^{(3)}))^\vee$.
    \end{enumerate}
\end{proposition}

\begin{proof}
    Set $\Delta^{[2]} := {\Delta}/{\Delta^{(2)}}$ and $V^{[3]}:=
    V/V^{(3)}$.  Because the presentation is minimal,
    $\Delta^{[2]}\simeq{V}/{V^{(2)}}$. We then have the following
    commutative diagram
    \begin{equation*}
        \xymatrix{1 \ar[r] & R \ar[r] \ar[d] & V \ar[r]
        \ar[d] & \Delta \ar[d] \ar[r] & 1 \\ 1 \ar[r] &
        \frac{V^{(2)}}{V^{(3)}} \ar[r] & V^{[3]} \ar[r] &
        \Delta^{[2]} \ar[r] & 1.}
    \end{equation*}
    Passing to the natural maps induced by the Hochschild-Serre-Lyndon
    spectral sequence with coefficients in $\Fp$, we obtain the
    commutative diagram
    \begin{equation*}\Small
        {\xymatrix{H^1(\Delta^{[2]}) \ar[r]^{\inf} \ar[d] &
        H^1(V^{[3]}) \ar[r]^-{\res} \ar[d] &
        H^1\left(\frac{V^{(2)}}{V^{(3)}}\right)^{\N^{[2]}}
        \ar[r]^-{\tra} \ar[d] & H^2(\Delta^{[2]}) \ar[r] \ar[d]
        & \dots \\ H^1(\N) \ar[r]^{\inf} & H^1(V)
        \ar[r]^{\res} & H^1(R)^\Delta \ar[r]^{\tra} &
        H^2(\Delta) \ar[r] & 0.}}
    \end{equation*}
    Since the extension $1 \to V^{(2)}/V^{(3)} \to V^{[3]} \to
    \N^{(2)} \to 1$ is central, $H^1(V^{(2)}/
    V^{(3)},\Fp)^{\Delta^{[2]}} = H^1(V^{(2)}/ V^{(3)}, \Fp)$.
    Additionally using the fact that the inflation map on the first
    cohomology group in each row of the previous diagram is an
    isomorphism, we may extract the following commutative square,
    with the rightmost transgression map an isomorphism:
    \begin{equation*}
        \xymatrix{H^1\left(\frac{V^{(2)}}{V^{(3)}},\Fp\right)
        \ar[d]^{\tra} \ar[r] & H^1(R,\Fp)^\N \ar[d]^{\tra}_{\simeq} \\
        H^2(\Delta^{[2]},\Fp) \ar[r] & H^2(\N,\Fp) }
    \end{equation*}
    Since $H^1(R,\Fp)^\N \simeq (R/R^p[R,V])^\vee$, we have
    \eqref{it:n1}.

    The leftmost transgression map, however, is also an isomorphism
    by \cite[1.1 and proof]{Ho}.  Now consider the natural map
    $R/R^p[R,V]\to V^{(2)}/V^{(3)}$ of abelian topological groups of
    exponent $p$.  The image of the natural map
    $H^1(V^{(2)}/V^{(3)},\Fp)\to H^1(R,\Fp)^{\N}$ is
    $H^1(RV^{(3)}/V^{(3)},\Fp)$. We then factor the horizontal maps
    of the commutative square into homomorphisms followed by
    injections:
    \begin{equation*}
        \xymatrix{H^1\left(\frac{V^{(2)}}{V^{(3)}},\Fp\right)
        \ar[d]^{\tra}_{\simeq} \ar[r] &
        H^1(\frac{RV^{(3)}}{V^{(3)}},\Fp) \ar@{^{(}->}[r] &
        H^1(R,\Fp)^\N \ar[d]^{\tra}_{\simeq} \\
        H^2(\Delta^{[2]},\Fp) \ar[r] & H^2(\N,\Fp)^\dec
        \ar@{^{(}->}[r] & H^2(\N,\Fp). }
    \end{equation*}
    We obtain isomorphisms
    \begin{equation*}
        H^1(R/(R\cap V^{(3)}),\Fp) \simeq H^1(RV^{(3)}/V^{(3)},\Fp)
        \simeq H^2(\N,\Fp)^\dec.
    \end{equation*}
    Hence $H^2(\N,\Fp)^\dec \simeq (R/(R\cap V^{(3)}))^\vee$, and
    we have proved \eqref{it:n2}.
\end{proof}

\begin{proposition}\label{pr:gammah}
    Let $\G$ and $H$ be pro-$p$-groups with maximal closed subgroups
    $\N$ and $N$, respectively, and $\alpha:\G \to H$ a surjection
    with $\alpha(\N) = N$ and $\Ker \alpha \subset \N^{(3)}$.  Write
    $G$ for $\G/ \N \simeq H/N$.  Then as $\Fp G$-modules,
    $H^2(\N,\Fp)^\dec \simeq H^2(N,\Fp)^\dec$.
\end{proposition}

\begin{proof}
    We prove that $H^2(\N,\Fp)^\dec \simeq H^2(N,\Fp)^\dec$ under a
    natural isomorphism, and it will follow that the isomorphism is
    $\Fp G$-equivariant.  Let $\theta: V \to \N$ be a minimal
    presentation of $\N$ with kernel $R$. We may choose a section
    $W\subset V^{(3)}$ of $\ker \alpha$ under the surjection
    $\theta$.  We obtain a minimal presentation of $N = \alpha(\N)$:
    $1 \to RW \to V \xrightarrow{\psi} N \to 1$, where $\psi =
    \alpha \circ \theta$.  By Proposition~\ref{pr:a2}\eqref{it:n2},
    $H^2(\N,\Fp)^\dec \simeq (RV^{(3)}/V^{(3)})^\vee \simeq
    (RWV^{(3)}/{V^{(3)}})^\vee \simeq H^2(N,\Fp)^\dec$.
\end{proof}

\begin{lemma}\label{le:lem}
    Let $p$ be prime and $\G$ a nonfree pro-$p$-group which is the
    absolute Galois group of a field $F$.  Then $\chr F\neq p$ and
    $\xi_p\in F$.
\end{lemma}

\begin{proof}
    Since $\G$ is not free, then by Witt's Theorem, $\chr F\neq p$.
    Since $([F(\xi_p):F],p)=1$, $\xi_p\in F$.
\end{proof}

\begin{proposition}\label{pr:deltapgroup}
    Let $p>3$ be prime.  Suppose that $\G$ is a pro-$p$-group and
    $\N$ is a maximal closed subgroup of $\G$. If the
    $\Fp(\G/\N)$-module $H^2(\N,\Fp)^\dec$ contains a cyclic summand
    of dimension $i$ with $3\le i<p$, then $\G$ is not an absolute
    Galois group. Moreover, if $\G$ contains a normal closed
    subgroup $\Lambda\subset \N$ with $\G/\Lambda\simeq \Z/p^2$,
    then we may take $2\le i<p$ in the same statement.
\end{proposition}

\begin{proof}
    Suppose that $\G=G_F$ for some field $F$.  Then $\N = G_E$ for
    some $E/F$ of degree $p$.  Since $H^2(\N,\Fp)\neq 0$, $G_E$ is
    not a free pro-$p$-group and therefore neither is $G_F$
    \cite[Corollary 3, \S I.4.2]{SGC}. Lemma~\ref{le:lem} gives
    $\chr F \ne p$ and $\xi_p\in F$. By \cite[Theorem~11.5]{MeSu},
    $H^2(\N,\Fp)$ is decomposable.  Therefore by \cite[Theorem
    1]{LMS}, $H^2(\N,\Fp)^\dec $ contains no cyclic
    $\Fp(\G/\N)$-summand of dimension $i$ with $3\le i<p$, a
    contradiction.  Moreover, if $\Lambda\subset N$ is a normal
    closed subgroup and $\G/\Lambda\simeq \Z/p^2\Z$, then by
    \cite{A}, $\xi_p\in N_{E/F}(E)$.  Letting $E=F(\root{p}\of{a})$
    for some $a\in F^\times$, we obtain in $H^2(\G,\Fp)$ that
    $(a).(\xi_p)=0$.  By \cite[Theorem 1]{LMS}, $H^2(\N,\Fp)^\dec $
    contains no cyclic $\Fp(\G/\N)$-summand of dimension $2$, again
    a contradiction.
\end{proof}

\begin{corollary}\label{co:deltan}
    Let $p>3$ be prime.  Suppose that $\G$ and $H$ are
    pro-$p$-groups with respective maximal subgroups $\N$ and $N$,
    and $\alpha:\G \to H$ a surjection with $\alpha(\N) = N$ and
    $\Ker \alpha \subset \N^{(3)}$.  Write $G$ for $\G/\N \simeq
    H/N$. If either $H^2 (\N,\Fp)^\dec$ or $H^2(N,\Fp)^\dec$
    contains a cyclic $\Fp G$-summand of dimension $i$ with $3\le
    i<p$, then neither $\N$ nor $H$ is an absolute Galois group.

    Moreover, suppose additionally that $H$ contains a normal closed
    subgroup $\Lambda\subset N$ with $H/\Lambda\simeq \Z/p^2\Z$.
    Then if either $H^2(\N,\Fp)^\dec$ or $H^2(N,\Fp)^\dec$ contains
    a cyclic $\Fp G$-summand of dimension $2$, then $H$ is not an
    absolute Galois group.
\end{corollary}

\begin{proof}
    By Proposition~\ref{pr:gammah} we have that $H^2(\N, \Fp)^\dec
    \simeq H^2(N,\Fp)^\dec$ as $\Fp G$-modules. The remainder
    follows from Proposition~\ref{pr:deltapgroup}.
\end{proof}

\begin{proof}[Proof of Theorem~\ref{th:fam}]
    Let $\Omega$ be the group of Proposition~\ref{pr:const}. Observe
    that $\N :=((\Omega\star\Sigma)\times p\Z_p)/\Ec$ is a maximal
    closed subgroup of $\G$ of index $p$.  Let $G=\G/\N$ and note
    that the actions of $G$ and $C$ on $\N$ are identical. By
    Corollary~\ref{co:deltan} it is enough to show that
    $H^2(\N,\Fp)^\dec$ for $\Ec=1$ contains a cyclic $\Fp G$-summand
    $M_i$ with $3\le i<p$. Assume then that $\Ec=1$. By
    \cite[Theorem~4.1.4]{NSW} we have $H^1(\N,\Fp) \simeq
    H^1(\Omega,\Fp) \oplus H^1(\Sigma,\Fp) \oplus H^1(p\Z_p,\Fp)$
    and $H^2(\N,\Fp) \simeq H^2(\Omega,\Fp) \oplus H^2(\Sigma,\Fp)
    \oplus H^1(\Omega,\Fp)\oplus H^1(\Sigma,\F_p)$. By
    Proposition~\ref{pr:const}, $H^2(\Omega,\Fp)$ is decomposable
    and so $H^2(\Omega,\Fp)$ is a direct summand of
    $H^2(\N,\Fp)^\dec$. From Proposition~\ref{pr:const}, we obtain
    that $H^2(\Omega,\Fp)$ contains an $\Fp G$-summand $M_i$ with
    $3\le i\le p-1$.
\end{proof}

\begin{remark*}
    The proof excludes the groups $\G$ from the class of absolute
    Galois groups by using the $n=2$ case of \cite[Theorem 1]{LMS}.
    However, neither a direct application of the $n=1$ case nor
    Theorem~\ref{th:t} excludes the groups $\G$. Observe also that
    the fact that $\Omega\rtimes \Z_p$ is not an absolute
    Galois group could be deduced from the main results in
    \cite{K}, using different methods. However, the fact that each
    $((\Omega\times \Sigma)\rtimes \Z_p)/\Ec$ is not an absolute
    Galois group does not follow from \cite{K}.
\end{remark*}

\section*{Acknowledgement}

We thank the referee of \cite{LMS} for suggesting that we illustrate
\cite[Theorem 1]{LMS} by constructing examples of pro-$p$-groups
which are not absolute Galois groups, for this appendix is the
result. The third author would also like to thank Sunil Chebolu,
Ajneet Dhillon and Vahid Shirbisheh for stimulating discussions
related to this paper.

\end{document}